# A Unique and Novel Optimization Model for Solving the Timetabling Problem at Finnish Universities

Kimmo Nurmi[a]*, Jari Kyngäs[a] and Pauliina Halme[b]

[a]*Satakunta University of Applied Sciences, Pori, Finland*

**cimmo.nurmi@samk.fi*

[b]*Eduix Ltd., Tampere, Finland*

Kimmo Nurmi

- Research Director at Satakunta University of Applied Sciences
- Adjunct Professor of computational intelligence, algorithms and optimization at University of Turku
- CEO, Chairman & Founder at CI Computational Intelligence Ltd

Jari Kyngäs

- Specialist Researcher at Satakunta University of Applied Sciences
- PhD in computer science, University of Turku
- CDO & Founder at CI Computational Intelligence Ltd

Pauliina Halme

- Partner Director & Partner at Eduix Ltd
- M.Sc. Admin., University of Tampere

# A Novel Optimization Model for Solving the Timetabling Problem at Finnish Universities

We study a timetabling problem faced by Finnish universities, where the challenge is to balance the preferences of three key stakeholders: the university, teachers, and students. Educational timetabling is an essential task that ensures that everything runs effortlessly and collaboratively. Its efficiency significantly impacts the smooth operation of the institution's educational programs and its financial success. As universities have evolved and competition among them has intensified, the features of the university timetabling problem have changed. This paper describes an ongoing optimization project with three Finnish universities and a leading provider of educational management systems for universities. The main goal of this cooperation is not merely to streamline the creation of timetables, but to enable an entirely new approach to manage the entire planning process. The project consists of three six-month phases: the feasibility study, the preparation phase, and the implementation. During the feasibility study, several research-related gaps and novel practical goals were identified. To address these gaps and goals, this paper presents a unique CRONOPT model for tackling three targets simultaneously: effective use of university resources, teacher workload and well-being, and smooth progression of students' studies. The optimization solution developed through the project integrates directly with the market-leading educational management system for Finnish universities.

## 1. The strategic and operational needs to be addressed

We started our timetabling research first thirty years ago [1] and then again twenty years ago [2, 3]. Unfortunately, the business possibilities back then were challenging, so we turned our interest to professional sports league scheduling [4] and workforce scheduling [5]. We have continued to work in these research and business areas since then (see e.g. [6] and [7]). The successful results in these areas led us to revisit the issue of timetabling problems. This was due to the encouragement and support we received from some Finnish universities.

Twenty years of business experience have taught us to approach a new problem first from a strategic perspective, then from an operational and practical perspective, and only finally by considering a scientific solution. It is understandable that researchers are first and mostly interested in finding the best possible scientific solutions. However, working this way makes it difficult to bridge the gap between theory and practice.

One of the few papers discussing the lack of attention on real-world implementation of optimization solutions in timetabling applications is the one by Oude Vrielink et al. [8]. They state that timetabling applications comprise only a limited amount of state-of-the-art research in timetabling solutions, leading to a gap of twenty years or more between state-of-the-art algorithms in research and the implemented algorithms in timetabling software. This is unfortunate, and results in unused opportunities for both researchers and universities.

To the best of our knowledge, this is the first paper describing the entire process, from the strategic definition of the timetabling problem, all the way to the novel model for balancing workload, increasing well-being and enabling smooth progression of students' studies, and to the implemented method to be adopted.

Optimal course schedules have a direct and multifaceted impact on the operations of a university, ranging from students' graduation timelines to the efficient use of resources. Their importance is particularly emphasized as combinations of on-campus, online, blended, and cross-disciplinary learning become more common. The most significant impacts include the following:

(1) Smooth progression of studies. The best placement of courses in the curriculum enables students to progress smoothly through their studies. Students maintain sufficient capacity for study and well-being without an unreasonable workload. Predictable schedules further help students to succeed in their studies.

(2) Quality of teaching and teacher well-being. A functional, balanced, fair, predictable, and transparent timetabling enables teaching staff to provide high-quality teaching and maintain their well-being without an unreasonable workload.

(3) Optimization of resources. Optimal planning of facilities reduces idle time in rooms, and the fragmentation of teachers' working time.

(4) Cost-effectiveness and graduation rates. Smooth study paths and high graduation rates directly impact university's funding, which in Finland is mainly based on degrees completed on time.

The paper is organized as follows. Section 2 describes the timetabling process in Finnish universities. In Section 3, we first review the most common school and university timetabling problems. Then we detail the specific characteristics of the problem faced by Finnish universities arising from the process described in the earlier section. Section 4 introduces a unique and novel model for balancing workload, increasing well-being and enabling smooth progression of students' studies that can be widely applied to almost any timetabling problem. In Section 5, we review the most popular methods used to solve timetabling problems. We then describe the PEASTP method which will be fine-tuned to solve the Finnish university timetabling problem. Finally, Section 6 describes the project's upcoming preparation and implementation phases.

## 2. Timetabling process

Work on developing the solution began in collaboration with three universities and a leading provider of educational management systems for Finnish universities. The main goal of this solution was not merely to streamline the creation of timetables, but to enable an entirely new approach to manage the entire planning process based on mutually agreed-upon principles and transparent quality criteria. In addition, an important goal for the participating universities is to set themselves apart from other universities and gain a competitive advantage.

The ongoing project consists of three six-month phases: the feasibility study, the preparation phase, and the implementation. The goal of the feasibility study was to design a service that would comprehensively overhaul the operational model for timetable planning. The universities defined the target state for optimization-assisted timetabling, preparation and planning, and the subsequent implementation of timetables. This target state combines modern optimization technology, an understanding of the universities' operating environment, and the system-wide integrations into a single, manageable whole. In this entirely new system,

- the university's strategic goals are considered
- planning is based on the university's pedagogical and operational guidelines
- the smooth progression of students' studies is a key objective
- teachers' workloads are systematically assessed from the perspective of well-being
- the solution integrates directly with the market-leading educational management system for Finnish universities
- optimization should be genuinely goal- and rule-based.

In the target state, timetables are generated using an optimization algorithm based on rules, objectives, and quality criteria defined by the university itself. This allows the university to differentiate itself from its peers through strategic priorities related to the implementation of educational activities and instructions.

The feasibility study revealed that the timetabling problem in Finnish universities, as defined mathematically, is a high school timetabling problem (see e.g. [1] and [2]). The problem can be defined in terms of the availability of teachers and rooms, a number of lessons to be taught by teachers to specific classes or students and a set of constraints. The construction of a timetable involves assigning resources, such as timeslots, teachers, student groups and rooms to a collection of lessons while minimising the constraint violations. We will give a description of the problem in Section 3.

***The process***

During the six-month feasibility study, we were able to outline a planning process that all three universities can apply, and we believe, all other Finnish universities as well. The main process consists of nine phases:

(1) Involve the university community in the principles and general guidelines for planning timetables, with the aim of
    (a) improving graduation rates,
    (b) optimizing resources,
    (c) ensuring transparency, fairness, and equality,
    (d) developing a model for balancing workload, increasing well-being and enabling smooth progression of students' studies, and
    (e) fostering a better overall understanding of the planning and optimization process from the perspective of all staff and students.

(2) Decisions by the university directors regarding the principles, general guidelines, and practices for timetable planning that guide timetabling planners.
(3) The choices made by those responsible for planning regarding the targets to be optimized, and their order of priority, as outlined in the model for balancing workload, increasing well-being and enabling smooth progression of students' studies.
(4) Decisions made by the heads responsible for the teaching process regarding the number and duration of course lectures, as well as their allocation to student groups, teachers, and rooms (recall here that we have a high-school timetabling problem in hand).
(5) Analyses conducted by timetabling planners using an analyser tool, which provide numerical data and traffic-light ratings indicating the quality of timetables that can be created based on the provided input data and given parameters.
(6) Reviewing the results of the analysis with the responsible parties and making the necessary changes to improve the expected quality of timetables.
(7) Optimization runs conducted by the timetabling planners until the results meet the target quality.
(8) Publication of the timetables.
(9) Audits and feedback on the timetables in accordance with the continuous improvement process + back to phase (1).

The first four phases involve collecting and selecting the correct and timely input data for optimization, and the next three phases involve the optimization itself. The fifth phase includes such analyses as whether there are mathematically sufficient (or excessive) facilities available, and whether teachers' requests are statistically feasible. The data for the first four phases are collected while updating, reviewing and approving curricula followed by annual course planning. Important decisions are made at various stages of the course planning. The most important decision regarding the optimization is the following:

> In Finland, an academic year is divided into 5−7 timeframes. Do we optimize all the weeks of the timeframe at the same time, even though the weeks most likely will end up having quite different outset?
> Or do we first optimize all the weeks the same, and then, in different weeks, we will move the necessary lectures from one time slot to another?

The first approach is able to consider the specific characteristics of each upcoming week, but for understandable reasons, it is mathematically more challenging to solve. Therefore, the second approach yields a better result in terms of qualitative criteria. In addition, the timetables would be more predictable for students because they are fairly static from week to week. As a result of the feasibility study, we decided to use the second approach.

Another important decision in course planning is to define the availability of teachers, student groups, and rooms. Some of the teachers also teach at other educational institutions, and some of them may live quite far from the campuses, so commuting times must also be considered. Similarly, the campus may also host lectures from other educational institutions, in which case some of the rooms are shared. For the student groups, particular attention should be paid to students' work placements at companies and other organizations.

One further issue to be decided is the determination of room groups. It would be computationally easier for the optimization if we could group similar rooms together.

The rooms in the same group are sufficiently equal to one another in terms of size, facilities, and suitability for teaching purposes.

## 3. The timetabling problem to be solved

Over the past decade, the number of publications on educational timetabling problems has clearly increased. Numerous variants of timetabling problems have been presented, proposed, and studied in the literature. They differ from one another depending on the type of educational institution in question (e.g. school, high-school, vocational school, and university), the type of meeting (e.g. classes, lessons, lectures, and exams) and the type of constraints, preferences, and optimization targets. Even the simplest forms of the problem types are NP-hard [9, 10].

The four main problem types are school timetabling, curriculum-based course timetabling, post-enrolment course timetabling, and exam timetabling. McCollum et al. [11] provide a good overview of the last three problem types. The exam timetabling problem involves assigning a specific set of exams into a limited number of timeslots and physical rooms while minimizing overlapping exam slots for popular course combinations, thereby maximizing the choices available to the student. From the perspective of this paper, we are interested in the other three problem types.

The pioneering work on course timetabling was the one by Hertz [12]. More recent good articles and reviews include [13, 14]. Curriculum-based course timetabling [15] creates a schedule before student enrolment based on predefined study programs, whereas Post-Enrolment timetabling [16] builds the schedule after students have registered to the courses, ensuring individual students have as few scheduling conflicts as possible.

In other words, Curriculum-based timetabling consists in scheduling courses to rooms and timeslots according to the curricula of the university. Post-Enrolment based timetabling is used when there is no concept for a curriculum, and the objective is to schedule courses into rooms and timeslots based on the selections made by the students. Thus, the main difference between these two models is that, in the curriculum model, a set of courses constitute the entire course load for a specific group of students, as in the post-enrolment model the restrictions and objectives are based on students' course enrolment.

With specific regard to school timetabling, the first seminal work was by de Werra [17]. Post et al. [18] summarized a set of common constraints for several school timetabling problems in different countries. Other notable papers include [10] and [19]. School timetabling problems can be defined in terms of (teacher, class, room) triples, and a number of lessons to be taught to these triples. The construction of a timetable involves assigning timeslots to the lessons while minimizing the constraint violations.

The primary difference between the two course timetabling problems and the school timetabling problem lies in student grouping and course enrolment. The students in school timetabling are bound to fixed student groups (or classes) that move together all day. In addition, school timetables are identical every week of the term, while university timetables can change weekly due to modular courses, intensive seminars, or bi-weekly labs.

***The Timetabling Problem in Finnish Universities***

Due to differences in cultures and educational systems, the challenges associated with university timetabling vary from country to country in terms of requirements and timetabling structures. We next study a timetabling problem faced by Finnish universities. The problem involves balancing the preferences of three key stakeholders, the university, teachers, and students.

In our implementation, a lecture can have multiple teachers, multiple student groups, and multiple rooms, or it can be taught online. It is obvious, that this makes timetabling particularly challenging. In fact, a lecture is not assigned to a specific room, but rather to a room group. A room group consists of rooms that are sufficiently similar in size and facilities, meaning that a lecture can be scheduled in any of the rooms belonging to that room group.

As stated in Section 2, the feasibility study revealed that the timetabling problem in Finnish universities, as defined mathematically, is a high school timetabling problem. The problem is to schedule predefined lectures of courses, consisting of (teachers, student groups, rooms), to the timeslots while balancing workload and maximizing teacher well-being and smooth progression of students' studies.

The strategic and process-oriented background of timetabling in Finnish universities were described in the first two sections. Next, we will examine some practical aspects of the operating environment and operations of the Finnish universities that must be considered when building the timetables. Finnish universities have multiple campuses, so transition times between buildings must be considered when scheduling lectures. In addition, some teachers live quite far from the campuses, so commuting times must also be considered.

There are many factors to consider when defining how different lectures of different courses are scheduled, with the aim of ensuring the smooth progression of students' studies. First, we need to determine the number of lectures of each course for each student group per week and the duration of each, for example twice a week for three hours at a time, or one two-hour lecture and one four-hour lecture. In addition, we need to define the constraints regarding the lectures for the same and different courses. These include lectures that cannot be scheduled on the same day, or on the next day.

Some lectures need to be in the correct order during the week. These are familiar issues that have been considered many times in the scheduling literature.

A new issue is to balance the cognitive load per workday and per study day. Some courses are more demanding and load-bearing than others. Another new issue is to prevent alternating between on-campus and online lectures on the same day. This will improve time management by reducing logistical stress, minimize cognitive fatigue, and prevent severe disruptions to student focus and engagement. The third issue is to find a good balance in the number of free periods (gaps) in students' study days. Some gaps are necessary so that students have access to special rooms (labs etc.) for individual studying during the day, and also to allow long enough lunch breaks. Teachers most likely prefer gapless timetables enabling them to have more uninterrupted time for research.

## 4. Unique and novel optimization model

University's timetables have a significant impact on the smooth operation of the institution's educational programs, as well as on its financial success. Competition for students has increased significantly because of economic challenges leading to pressure to cut back on spaces, facilities, and other resources. The efficient use of campuses and their spaces is one of the most significant factors affecting university costs. This is also related to the cost-reducing opportunities offered by distance learning.

In addition, well-structured timetables have a direct impact to the satisfaction and well-being of students and staff. This section presents a unique model for tackling three targets simultaneously: effective use of university resources, teacher workload and well-being, and smooth progression of students' studies. Only a few papers have considered quality issues of timetables. Schimmelpfeng and Helber [20] described the implementation process and reported results from an anonymous satisfaction survey with respect to their planning approach. Fairness objective in addition to the penalty objective was studied in [21].

***Learning from the WTTL model***

The timetabling model presented in this section derives from the Working Time Traffic Light (WTTL) model [22] by the Finnish Institute of Occupational Health (FIOH), which operates under the Ministry of Social Affairs and Health. They have more than 20 years of experience in studying, examining, and assessing the effects and consequences of shift work. Their recommendations are based on the daily working time statistics of more than 300 000 employees working in healthcare and social services sector.

Our research and business operations in workforce scheduling are closely related to the FIOH studies. Our staff rostering optimization is based on their model, and it has been integrated into a market-leading workforce management software in Finland (see e.g. [5]). Similarly, the optimization solution developed as a result of the ongoing work described in this paper integrates directly with the market-leading educational management system for Finnish universities.

The timetabling model is based entirely on the research findings underlying the WTTL model. The most important ones are the following:

(1) Well-planned schedules can reduce work-related stress and promote recovery from work.
(2) Good work schedules support work ability and long careers.
(3) Work schedules should respect peoples' internal clocks.
(4) People should have sufficient control over their working hours to ensure good work-life balance.

Stress factors include the predictability of work, total working hours, overtime, the length of the workday, the number of evening hours, daily and weekly recovery time and work breaks. Short rest time between two working days should be avoided, as ruled by the EU Working Time Directive. Quick returns have detrimental effects on acute health problems [23]. They encumber many detriments, such as short sleep duration, slightly prolonged sleep onset latency, more abrupt ending of main sleep period, increased sleepiness, and higher level of perceived stress on the following day [24].

Some of us are naturally early birds, while others are night owls. Thus, universities can increase productivity, as well as job satisfaction by respecting the teachers' body clocks with better targeted timetables. Social activities occur mostly in the evenings and on weekends forming a social rhythm [25]. Working hours in the evenings impair social well-being and work-life balance [26].

Increase of working time control is a potential way to improve health, well-being and work participation [27]. High levels of work time control have been shown to positively influence mental health outcomes such as affective well-being and perceived stress [28]. The ability to influence the academic week, teaching days, the rhythm of lectures, and teaching rooms through personal preferences is also an advantage in the job market. However, it should be noted that teachers' preferences are accommodated as much as possible, provided they are mathematically feasible and do not compromise the overall quality of the timetable, the quality of the timetables for their student groups, or the quality of other teachers' timetables.

***The CRONOPT model***

We next present the unique and novel model for balancing workload, increasing well-being and enabling smooth progression of students' studies. We will call the model CRONOPT from now on.

From teacher perspective, work is associated with pressure on teaching quality, completion of students' degrees within the target timeframe, social and family life, motivation, and health issues. From student perspective, studying is associated with pressure on study progression, motivation, working alongside studies, social life, and health issues. As a result, the CRONOPT model should address these issues as part of balancing workload, increasing well-being and enabling smooth progression of students' studies. The timetabling problem and CRONOPT model is stated as follows:

**Timetable**

The initial state for optimization consists of number of lectures. Each lecture is assigned a course, the duration of the lecture, the teacher who teaches it, the student group it is aimed at, and the room in which it is lectured. Multiple number of teachers, student groups, and rooms can be assigned to a single lecture, although there is usually only one of each. In addition, recall the earlier discussion about room groups.

**Problem**

The timetabling problem is to assign a timeslot for each lecture while considering aspects of the three components: 1) the necessary conditions, 2) the optimization targets for teacher workload and well-being, and 3) the optimization targets for smooth progression of students' studies. Some lectures may be pre-assigned, and they timing cannot be rescheduled.

**Necessary conditions**

(N01) All lectures should be scheduled. Some lectures may be pre-assigned

(N02) Lectures must not overlap in time in terms of courses, teachers, student groups, and rooms.

(N03) Blocked timeslots for courses, teachers, student groups, and rooms must be adhered to.

**Optimization targets for teacher workload and well-being**

Timing of working hours

(TT01) The maximum number of hours per workweek should not be exceeded

(TT02) The maximum workload of consecutive workdays should not be exceeded

(TT03) The length of the shortest workday should not go below the minimum

(TT04) The length of the longest workday should not go above the maximum

(TT05) The longest number of consecutive hours should not go above the maximum

Recovery and rest times

(TR06) The number of evening hours per week should not go above the maximum

(TR07) The number of weekend hours per week should not go above the maximum

(TR08) The longest work period without a day off should not go above the maximum

(TR09) The daily rest in hours should not go below the minimum

(TR10) The next workday after working on Saturday should not be Monday

(TR11) The cognitive load (demanding and load bearing courses) per day should not go above the maximum

Scheduling logic and flow

(TS12) Lectures should be scheduled in the correct order

(TS13) On-campus and online courses should not be scheduled on the same day

(TS14) Courses that require time-consuming transitions from one building to another may should not be scheduled on the same day

(TS15) Lunch breaks should be scheduled according to the first teaching timeslot each day

Control over working hours

(TC16) Free periods (gaps) should, or should not, be scheduled according to personal preferences

(TC17) Personal internal body clocks (early bird / night owl) should be considered

(TC18) Personal priority preferences should be fulfilled

(TC19) Secondary preferences should be considered

**Optimization targets for smooth progression of students' studies**

Timing of study hours

(ST01) The maximum number of hours per study week should not be exceeded

(ST02) The maximum workload of consecutive study days should not be exceeded

(ST03) The length of the shortest study day should not go below the minimum

(ST04) The length of the longest study day should not go above the maximum

(ST05) The longest number of consecutive hours should not go above the maximum

Recovery and rest times

(SR06) The number of evening hours per week should not go above the maximum

(SR07) The number of weekend hours per week should not go above the maximum

(SR08) The longest work period without a day off should not go above the maximum

(SR09) The daily rest in hours should not go below the minimum

(SR10) Free periods (gaps) should, or should not, be scheduled according to student group preferences

(SR11) The cognitive load (demanding and load bearing courses) per day should not go above the maximum

Scheduling logic and flow

(SS12) Lectures should be scheduled in the correct order

(SS13) Some pairs of lectures should not be scheduled on the same day

(SS14) Some pairs of lectures should not be scheduled on consecutive days

(SS15) On-campus and online courses should not be scheduled on the same day

(SS16) Courses that require time-consuming transitions from one building to another may should not be scheduled on the same day

(SS17) Lunch breaks should be scheduled according to the first teaching timeslot each day

(SS18) Access to special rooms (labs etc.) for individual studying

These 40 constraints and optimization targets involve a vast number of interdependent, indirect, and hidden relationships that exist between different planning objectives. It is not mathematically possible to achieve all objectives, but through careful planning, analysis, and the proper use of optimization, it is possible to create timetables that best fulfil the wishes of all three parties: university, teachers, and students.

For teachers, optimization enables effective, balanced, fair, predictable, and transparent schedules. This ensures high-quality teaching and the well-being of teaching staff without an unreasonable workload. However, we must understand that teachers may have a different view of what fairness or equality means for them.

For students, optimization enables smooth progress in their studies. Students can maintain their ability to study, as well as their well-being without an unreasonable study load. A motivating student experience ensures that degrees are completed in time and with a high graduation rate.

Following the WTTL model, the CRONOPT model uses four levels in to define recommended values for optimization targets. The levels are safe load, increased overload, overload, and heavy overload. It is the responsibility of the university community to define the threshold values of each optimization target in collaboration with the HR department, the Joint Consultative Committee (JCC), staff, and student organizations. An example of the implementations of the CRONOPT model for some optimization targets is given in Table 1.

| **CRONOPT MODEL** | | | | | |
|---|---|---|---|---|---|
| **Teacher workload and well-being** | | **Safe** | **Increased** | **Overload** | **Heavy** |
| TT04 | Length of the longest workday | ≤ 6 | 7 | 8 | ≥ 9 |
| TR09 | Daily rest in hours | ≥ 14 | 13 | 12 | ≤ 11 |
| TS13 | Campus-to-online daily transitions | 0 | ≤ 2 % | ≤ 5 % | > 5 % |
| TC18 | Priority preferences to be fulfilled (personal) | > 90 % | 75-90 % | 50-75 % | < 50 % |
| **Smooth progression of students' studies** | | **Safe** | **Increased** | **Overload** | **Heavy** |
| ST05 | Longest number of consecutive hours | ≤ 4 | 5 | 6 | ≥ 7 |
| SR11 | Cognitive load per day | ≤ 10 | 11-12 | 13-14 | ≥ 15 |
| SS13 | Number of unsuitable lectures on the same day | 0 | ≤ 2 % | ≤ 5 % | > 5 % |
| SS18 | Access to special rooms for individual studying | No viol. | ≤ 2 % viol. | ≤ 5 % viol. | > 5 % viol. |

Table 1. An example of the CRONOPT model implementation for some optimization targets.

## 5. Solution method

The number of publications on solution methods for timetabling problems has clearly increased in recent years. Good reviews of the popular methods include [10], [29] and [30]. The solution methods can be classified primarily into local search methods, population-based metaheuristics, hyper-heuristics, and hybrid methods.

Local search methods (single solution-based meta-heuristics) [49] can address large-scale problem instances in reasonable computational times. The three most celebrated local search methods are tabu search [32] and simulated annealing [33], and variable neighbourhood search [34]. They are quite simple to implement, and still effective in achieving high quality solutions. However, researchers need to consider parameter tuning when choosing local search approaches. Researchers are working on designing an optimization algorithm that is not only effective but requires less manual parameter setting.

Population-based meta-heuristics operate on a population of solutions and apply various operators and rules to evolve a new population of solutions in the neighbourhood areas of current ones. Examples of population-based meta-heuristics are genetic algorithms [35], ant colony optimization [36] and particle swarm optimization [37]. Examples of the use of these methods in timetabling can be found in [38], [39] and [40]. Population-based methods are superior compared to others in terms of solution space exploration. However, one of the drawbacks of these approaches are the computational times required in finding good quality solutions.

Due to the limitation of meta-heuristic approaches which require intensive parameter tuning, hyper-heuristic approaches were introduced [41]. Hyper-heuristics are general, simple, and fast algorithms applicable to variety of problem domain and can event adapt to different instances of a given benchmark dataset. Hyper-heuristics are heuristics to choose heuristics (algorithms), working on a search space of heuristics (algorithms) instead of a search space of solutions. The challenges in hyper-heuristic approaches are balancing information exchange and maintaining a problem domain barrier between the low-level heuristics and the high-level search methodology.

***PEASTP***

The PEASTP optimization algorithm [42] will be fine-tuned to solve the Finnish university timetabling problem. It has been in successful commercial use for over 20 years in sectors such as social and health care, public transportation, call centers, and professional sports (see e.g. [6] and [7]). The effectiveness of the optimization has also been demonstrated in publicly available scientific performance tests, and over 70 international research articles have been published on the algorithm’s performance. However, it should be clear that even though an algorithm can be powerful for several problem types, it may perform poorly for other problem types. The No Free Lunch Theorem [43] implies that there cannot exist a superior optimization method.

The PEASTP algorithm uses a population of solutions in each iteration. The worst solution is replaced with the best one at given intervals, i.e. elitism is used. The ejection chain local search explores promising areas in the search space. The ejection chain search is improved by introducing a tabu list, which prevents reverse order moves in the same sequence of moves. A simulated annealing refinement is used to decide whether to commit to a sequence of moves in the ejection chain search.

The simulated annealing refinement and tabu search are used to avoid staying stuck in promising search areas too long. Shuffling operators assist in escaping from local optima. They are used to perturb a solution into a potentially worse solution in order to escape from local optima. A shuffling followed by several ejection chain searches obtains better solutions using the same idea as the ruin and recreate method. The heuristic uses a penalty method, which assigns dynamic weights to the hard constraints based on the constant weights assigned to the soft constraints.

## 6. Preparation and Implementation Phases

Many publications discuss the complete automation of university course schedules. However, human input is still needed to guide and validate the search process. Twenty years of business experience have taught us to approach a new problem first from a strategic perspective, then from an operational and practical perspective, and only finally by considering a scientific solution. The goal in getting a new optimization solution into production is not merely succeeding in optimization algorithm itself, but to enable an entirely new approach to manage the entire planning process based on mutually agreed-upon principles, and transparent quality criteria.

In this paper, we described a project for solving the Finnish university timetabling problem consists of three six-month phases: the feasibility study, the preparation phase, and the implementation. and students. In the feasibility study, we have identified several research-related gaps and novel practical goals. We presented a unique CRONOPT model for tackling three targets simultaneously: effective use of university resources, teacher workload and well-being, and smooth progression of students' studies.

In the ongoing preparation phase, the PEASTP optimization algorithm will be fine-tuned to solve the Finnish problem. We are currently developing a new instance generator for versatile school timetabling problems. Instances are generated according to the CRONOPT model, and the efficiency of the PEASTP algorithm will be evaluated based on how well it solves these instances. The generator is not designed for algorithm competitions; rather, it is a tool for practical-oriented researchers who want to develop and improve their own algorithms to solve practical school timetabling instances. For practitioners, the generator can also be used to find hidden connections between optimization targets and in their mutual order of importance.

In the implementation phase, the CRONOPT approach and the PEASTP algorithm are utilized to solve real-world cases obtained from partner universities. We believe that the computational results will reveal clear trade-offs as well as benefits between the main optimization goals and withing the single optimization targets. It would not be surprising to find out that improvements in one target within one main goal comes at the expense of some others. The reason is the multifactorial and interdependent relations of the optimization targets, i.e. the targets compete with each other.

The optimization solution developed as a result of the project integrates directly with the market-leading educational management system for Finnish universities. We believe that the new planning process, the CRONOPT model, and the PEASTP method is an efficient trio for generating timetables for Finnish universities. The integrated solution will provide timetabling planners a comprehensive tool for implementing the goals and guidelines jointly and transparently agreed upon by the university community, analyze the feasibility of their implementation, compare the effects of priorities regarding different optimization targets, and ultimately automate the creation of the timetables.

**Conflicts of interest**

The authors declare that there are no conflicts of interest regarding the publication of this paper.

**Generative AI statement**

Authors declare that no generative AI was used in creating content for this paper, including writing, data analysis, or creation of figures. The paper is the original product of human intellectual effort, and the work conforms to traditional standards of human authorship and academic integrity.